\input amstex
\magnification=\magstep1 
\baselineskip=13pt
\documentstyle{amsppt}
\vsize=8.7truein \CenteredTagsOnSplits \NoRunningHeads
\def\zz{\bold{z}}

\def\tr{\operatorname{trace}}
\def\rk{\operatorname{rank}}

\def\ind{\operatorname{ind}}

\def\NN{\Cal{N}}
\def\ind{\operatorname{Ind}}
\topmatter
 
\title Integrating products of quadratic forms   \endtitle 
\author Alexander Barvinok  \endauthor
\address Department of Mathematics, University of Michigan, Ann Arbor,
MI 48109-1043, USA \endaddress
\email barvinok$\@$umich.edu \endemail
\date February  2020 \enddate
\thanks  This research was partially supported by NSF Grant DMS 1855428. 
\endthanks 
\keywords quadratic equations, algorithm, interpolation method, integration \endkeywords
\abstract 
We prove that if $q_1, \ldots, q_m: {\Bbb R}^n \longrightarrow {\Bbb R}$ are quadratic forms in variables $x_1, \ldots, x_n$ such that  each $q_k$ depends on at most $r$ variables 
and each $q_k$ has common variables with at most $r$ other forms, then the average value of the product $\left(1+ q_1\right) \cdots \left(1+q_m\right)$ with respect to the standard Gaussian measure in ${\Bbb R}^n$ can be approximated within relative error $\epsilon >0$ in quasi-polynomial $n^{O(1)} m^{O(\ln m -\ln \epsilon)}$ time, 
provided $|q_k(x)| \leq \gamma \|x\|^2 /r$ for some absolute constant $\gamma > 0$ and $k=1, \ldots, m$. When $q_k$ are interpreted as pairwise squared distances for configurations of points in Euclidean space, the average can be interpreted as the partition function of systems of particles with mollified logarithmic potentials. We sketch a possible application  to testing the feasibility of systems of real quadratic equations.
 \endabstract
\subjclass 14Q30, 65H14, 68Q25, 68W25, 90C23 \endsubjclass
\endtopmatter
\document

\head 1. Introduction and main results \endhead

Integration of high degree multivariate polynomials is computationally difficult and no efficient algorithms are known except in few special cases, when the polynomials have a rather simple algebraic structure (close to a power of a linear form), cf. \cite{B+11}, or have some very nice analytic properties (slowly varying or, most notably, log-concave), cf. \cite{LV07}. Since a general $n$-variate polynomial $p$ of degree $d$ is defined by ${n +d \choose d}$ parameters (for example, coefficients), the problem becomes interesting for large $n$ and $d$ only if $p$ has some special structure (such as the product of low-degree polynomials), which allows us to define $p$ using much fewer parameters.

In this paper, we integrate products of quadratic forms with respect to the Gaussian measure in ${\Bbb R}^n$.
We relate the problem to partition functions of mollified logarithmic potentials and to testing the feasibility of systems of real quadratic equations.

Our algorithms are deterministic and based on the method of polynomial interpolation, which has been recently applied to a variety of partition functions in combinatorial (discrete) problems, cf.
\cite{Ba16}. In continuous setting, the method was applied to computing partition functions arising in quantum models \cite{B+19}, \cite{H+19}.

\subhead (1.1) Quadratic forms on ${\Bbb R}^n$ \endsubhead We consider Euclidean space ${\Bbb R}^n$ endowed with the standard inner product 
$$\langle x, y \rangle = x_1 y_1 + \ldots + x_n y_n \quad \text{for} \quad x=\left(x_1, \ldots, x_n \right) \quad \text{and} \quad y=\left(y_1, \ldots, y_n \right)$$
 and corresponding Euclidean norm 
$$\|x\|=\sqrt{\langle x, x \rangle}=\sqrt{x_1^2 + \ldots + x_n^2} \quad \text{for} \quad x=\left(x_1, \ldots, x_n\right).$$
Let $q_1, \ldots, q_m: {\Bbb R}^n \longrightarrow {\Bbb R}$ be quadratic forms defined by
$$q_k(x)={1 \over 2} \langle Q_k x, x \rangle \quad \text{for} \quad k=1, \ldots, m, \tag1.1.1$$
where $Q_1, \ldots, Q_m$ are $n \times n$ real symmetric matrices. 

Our first result concerns computing the integral
$${1 \over (2 \pi)^{n/2}} \int_{{\Bbb R}^n} \left(1+ q_1(x)\right) \cdots \left(1+q_m(x)\right) e^{-\|x\|^2/2} \ dx. \tag1.1.2$$
The idea of the interpolation method is to consider (1.1.2) as a one-parameter perturbation a much simpler integral, in our case, of 
$${1 \over (2 \pi)^{n/2}} \int_{{\Bbb R}^{n}} e^{-\|x\|^2/2}\ dx =1. \tag1.1.3 $$
For the method to work, one should show that there are no zeros in the vicinity of a path in the complex plane which connects (1.1.2) and (1.1.3). We prove the following result.
\proclaim{(1.2) Theorem} There is an absolute constant $\gamma >0$ (one can choose $\gamma=0.151$) such that the following holds.
Let $q_k: {\Bbb R}^n \longrightarrow {\Bbb R}$, $k=1, \ldots, m$, be quadratic forms. Then 
$${1 \over (2 \pi)^{n/2}} \int_{{\Bbb R}^n} \left(1 + \omega q_1(x)\right) \cdots \left(1+\omega q_m(x)\right) e^{-\|x\|^2/2} \ dx \ne 0$$
for all $\omega \in {\Bbb C}$ such that $|\omega| \leq \gamma$, provided
$$\left| q_k(x) \right| \ \leq \ {1 \over \max\{m, n\}} \|x\|^2 \quad \text{for} \quad k=1, \ldots, m.$$
\endproclaim

By interpolation, for any constant  $0 < \gamma' < \gamma$, fixed in advance,
we obtain an algorithm which, given quadratic forms $q_1, \ldots, q_m: {\Bbb R}^n \longrightarrow {\Bbb R}$, computes (1.1.2) within relative error $0 < \epsilon < 1$ in quasi-polynomial 
$n^{O(1)} m^{O(\ln m -\ln \epsilon)}$ time provided
$$\left| q_k(x) \right| \ \leq \ {\gamma' \over \max\{m, n\}}\|x\|^2 \quad \text{for} \quad k=1, \ldots, m. \tag1.2.1$$
Note that by Theorem 1.2 and (1.1.3), the value of (1.1.2) is positive, as long as (1.2.1) holds.

Some remarks are in order. 

First, we note that the integrand in (1.1.2) can vary wildly. Indeed, for large $n$ the bulk of the standard Gaussian measure in ${\Bbb R}^n$ is concentrated in the vicinity of the sphere $\|x\|=\sqrt{n}$, see for example, Section V.5 of \cite{Ba02}. Assuming that $m=n$, we can choose $q_k(x) \sim \|x\|^2/n$ so that (1.2.1) is satisfied. Then, in the vicinity of the sphere $\|x\|=\sqrt{n}$, the product $\left(1 + q_1(x)\right) \cdots \left(1+q_m(x)\right)$ in (1.1.2) varies within an exponential in $m$ factor, and is not at all well-concentrated.

Second, if the quadratic forms $q_1, \ldots, q_m$ exhibit simpler combinatorics, we can improve the bounds accordingly. We prove the following result.
\proclaim{(1.3) Theorem} There is an absolute constant $\gamma > 0$ (one can choose $\gamma =0.151$) such that the following holds.
Let $q_k: {\Bbb R}^n \longrightarrow {\Bbb R}$, $k=1, \ldots, m$, be quadratic forms. Suppose further that each form depends on not more than $r$ variables 
among $x_1, \ldots, x_n$ and that each form has common variables with not more than $r$ other forms. Then 
$${1 \over (2 \pi)^{n/2}} \int_{{\Bbb R}^n} \left(1 + \omega q_1(x)\right) \cdots \left(1+\omega q_m(x)\right) e^{-\|x\|^2/2} \ dx \ne 0$$
for all $w \in {\Bbb C}$ such that $|\omega| \leq \gamma$, provided
$$\left| q_k(x) \right| \ \leq \  { \|x\|^2 \over r} \quad \text{for} \quad k=1, \ldots, m.$$
\endproclaim
By interpolation, for any constant  $0 < \gamma' < \gamma$, fixed in advance,
we obtain an algorithm which, given quadratic forms $q_1, \ldots, q_m: {\Bbb R}^n \longrightarrow {\Bbb R}$ as in Theorem 1.3 computes (1.1.2) within relative error $0 < \epsilon < 1$ in quasi-polynomial 
$n^{O(1)} m^{\ln m - \ln \epsilon}$ time provided
$$\left| q_k(x) \right| \ \leq \ {\gamma' \over r}\|x\|^2 \quad \text{for} \quad k=1, \ldots, m. \tag1.3.1$$
We prove Theorems 1.2 and 1.3 in Section 3 and describe the algorithm for computing (1.1.2) in Section 4.
In Section 2, we discuss connections with systems of particles with mollified logarithmic potentials and possible applications to testing the feasibility of systems of multivariate real quadratic equations.

\head 2. Connections and possible applications \endhead

\subhead (2.1) Partition functions of mollified logarithmic potentials \endsubhead Let $n=ds$ and let us interpret ${\Bbb R}^n ={\Bbb R}^d \oplus \cdots \oplus {\Bbb R}^d$ as the space of all
ordered $s$-tuples $\left(v_1, \ldots, v_s\right)$ of points $v_i \in {\Bbb R}^d$. Hence the distance between $v_i$ and $v_j$ is $\|v_i -v_j\|$. 

Let us fix some set $E$ of $m$ pairs $\{i, j\}$ of indices $1 \leq i < j \leq s$ and suppose that the energy of a set of points $\left(v_1, \ldots, v_s\right)$ is defined by 
$$-\sum_{\{i, j\} \in E} \ln \left(1 + \alpha \|v_i - v_j\|^2 \right) + {1 \over 2} \sum_{i=1}^n \|v_i\|^2, \tag2.1.1$$
where $\alpha >0$ is a parameter. The first sum in (2.1.1) indicates that there a repulsive force between any pair $\{v_i, v_j\}$ with $\{i, j \} \in E$ (so that the energy decreases if the distance between 
$v_i$ and $v_j$ increases), while the second sum indicates that there is a force pushing the points towards 0 (so that the energy decreases when each $v_i$ approaches $0$).
When $\alpha =0$, the repulsive force disappears altogether, and $\alpha \longrightarrow +\infty$, the repulsive force behaves as a Coulomb's force with logarithmic potential,
since
$$\lim_{\alpha \longrightarrow +\infty} \ln \left(1 + \alpha \|v_i - v_j \|^2\right) - \ln \alpha =2\ln \|v_i - v_j\|.$$
Thus the integral 
$${1 \over (2 \pi)^{n/2}} \int_{{\Bbb R}^n} \prod_{\{i, j\} \in E} \left(1+ \alpha \|v_i - v_j\|^2\right) e^{-(\|v_1\|^2 + \ldots + \|v_s\|^2)/2}  dx,  \tag2.1.2$$
which is a particular case of (1.1.2),
can be interpreted as the partition function of points with ``mollified" or ``damped" logarithmic potentials. One can think of (2.1.2) as the partition function for particles with genuine logarithmic potentials, provided each particle is confined to its own copy of ${\Bbb R}^d$ among a family of parallel $d$-dimensional affine subspaces in some higher-dimensional Euclidean space.

The integral (2.1.2) can be considered as a ramification of classical Selberg-type integrals
for logarithmic potentials:
$$\aligned &{1 \over (2 \pi)^{n/2}} \int_{{\Bbb R}^n} \prod_{1 \leq i < j \leq n} \left| x_i - x_j\right|^{2 \gamma} e^{-(x_1^2 + \ldots + x_n^2)/2} \ dx_1 \cdots d x_n \\&\quad = 
\prod_{j=1}^n {\Gamma(1+j\gamma) \over \Gamma(1+\gamma)}, \endaligned \tag2.1.3$$
see for example, Chapter 17 of \cite{Me04}. The integral (2.1.3) corresponds to points in ${\Bbb R}^1$ and a similar integral is computed explicitly for points in ${\Bbb R}^2$ (and $\gamma=1$), see Section 17.11 of \cite{Me04}. For higher dimensions $d$ no explicit formulas appear to be known.

In contrast, we compute integrals (2.1.2) approximately for certain values of $\alpha$, but we allow arbitrary dimensions and can choose an arbitrary set of pairs of interacting points (and we can even choose different $\alpha$s for different pairs of points). Theorem 1.3 can be interpreted as the absence of phase transition in the Lee - Yang sense \cite{YL52}, if $\alpha$ is sufficiently small. For example, if the set $E$ consists of all ${s \choose 2}$ pairs 
$\{i, j\}$, Theorem 1.3 implies that there is no phase transition (and the integral can be efficiently approximated) if 
$$\alpha \ < \ {\beta \over \max\{d, s\}}$$
for some absolute constant $\beta >0$.

\subhead (2.2) Applications to systems of quadratic equations \endsubhead Every system of real polynomial equations can be reduced to a system of quadratic equations, as one can successively reduce the degree by introducing new variables via substitutions of the type $z:=xy$. A system of quadratic equations can be solved in polynomial time when the number of equations is fixed in advance, \cite{Ba93}, \cite{GP05}, but as the number of equations grows, the problem becomes computationally hard.

 Here we are interested in the systems of equations of the type 
$$q_k(x) =1 \quad \text{for} \quad k=1, \ldots, m, \tag2.2.1$$
where $q_k: {\Bbb R}^n \longrightarrow {\Bbb R}$ are positive semidefinite quadratic forms. 
Such systems naturally arise in problems of distance geometry, where we are interested to find out if there are configurations of points in ${\Bbb R}^d$ with prescribed distances between some pairs of points and in which case $q_k$ are scaled squared distances between points, see \cite{CH88}, \cite{L+14} and Section 2.1. Besides, finding if a system of homogeneous quadratic equations has a non-trivial solution
$$q_k(x) = 0 \quad \text{for} \quad k=1, \ldots, m \quad \text{and} \quad \|x\|=1 \tag2.2.2$$
can be reduced to (2.2.1) with positive definite forms $q_k$ by adding $\|x\|^2$ to the appropriately scaled equations in (2.2.2).

Suppose that 
$$\sum_{k=1}^m q_k(x) = {\|x\|^2 \over 2}  \tag2.2.3$$
in (2.2.1). By itself, the condition (2.2.3) is not particularly restrictive: if the sum of in the left hand side of (2.2.3) is positive definite, it can be brought to the right hand side by an invertible linear transformation of $x$.

Let us choose an $\alpha > 0$ such that the scaled forms $\alpha q_k$ satisfy (1.3.1), so that the integral 
$${1 \over (2\pi)^{n/2}} \int_{{\Bbb R}^n} \left(1 + \alpha q_1(x)\right) \cdots \left(1 +\alpha q_m(x)\right) e^{-\|x\|^2/2} \ dx \tag2.2.4$$
can be efficiently approximated. We would like to argue that the value of the integral (2.2.4) can provide a reasonable certificate which allows one to distinguish systems (2.2.1) with many ``near solutions" $x$ from the systems that are far from having a solution.

We observe that the system (2.2.1) has a solution if and only if the system 
$$q_k(x)=t \quad \text{for} \quad k=1, \ldots, m \tag2.2.5$$
has a solution $x \in {\Bbb R}^n$ for any $t > 0$.

Let us find $0 < \beta < 1$ such that
$$2m\left({1 \over \beta} - {1 \over \alpha}\right)= {n \over 1-\beta}. \tag2.2.6$$
Indeed (2.2.6) always has a (necessarily unique) solution $0 < \beta < 1$, since for $\beta \approx 0$ the right hand side is bigger than the left hand side, while for $\beta \approx 1$ the left hand side is bigger than the right hand side. 

Because of (2.2.3), we can rewrite (2.2.4) as 
$${1 \over (2 \pi)^{n/2}} \int_{{\Bbb R}^n} e^{-{(1-\beta) \|x\|^2 \over 2}} \prod_{k=1}^m \bigl( 1+ \alpha q_k(x)\bigr) e^{-\beta q_k(x)} \ dx. \tag2.2.7$$
We observe that if $\alpha > \beta$ then the maximum value of 
$$\left(1 + \alpha t \right) e^{-\beta t} \quad \text{for} \quad t \geq 0$$
is attained at 
$$t={1 \over \beta} -{1 \over \alpha} \ > \ 0 \tag2.2.8$$ and is equal to 
$${\alpha \over \beta} \exp\left\{ {\beta \over \alpha} -1 \right\} \ > \ 1$$
and hence the maximum value of the product of the $m$ factors in (2.2.7) is 
$$\left({\alpha \over \beta} \exp\left\{ {\beta \over \alpha} -1 \right\}\right)^m$$
and attained if and only if the system (2.2.1) and hence (2.2.6) has a solution $x$.

Also, if $x$ is a solution to (2.2.5), by (2.2.3), (2.2.6) and (2.2.8), we have 
$$\|x\|^2=2tm={n \over 1-\beta}.$$
The Gaussian probability measure in ${\Bbb R}^n$ with density
$${(1-\beta)^{n/2} \over (2\pi)^{n/2}} e^{-{(1-\beta)\|x\|^2 \over 2}},$$
is concentrated in the vicinity of the sphere $\|x\|^2=n/(1-\beta)$,
cf., for example, Section V.5 of \cite{Ba05} for some estimates. Therefore, if for the system (2.2.1) there are sufficiently many ``near solutions" $x$, we should have the value of the integral (2.2.4) sufficiently close to 
$$\left( {\alpha \over \beta} \exp\left\{ {\beta \over \alpha} -1 \right\}\right)^m (1-\beta)^{-{n\over 2}},$$
while if the system (2.2.1) is far from having a solution, the value of the integral will be essentially smaller.

\head 2. Proofs of Theorems 1.2 and 1.3 \endhead

Choosing $r=\max\{m, n\}$, we obtain Theorem 1.2 as a particular case of Theorem 1.3. Hence we prove Theorem 1.3 only.

For a real symmetric $n \times n$ matrix $Q$ we denote 
$$\|Q\|=\max_{\|x\|=1} \| Q x \|$$
its operator norm.

We start with a simple formula, cf. also \cite{Ba93}. 

\proclaim{(3.1) Lemma} 
Let $q_1, \ldots, q_m: {\Bbb R}^n \longrightarrow {\Bbb R}$ be quadratic forms,
$$q_k(x)={1 \over 2} \langle Q_k x, x \rangle \quad \text{for} \quad k=1, \ldots, m,$$
where $Q_1, \ldots, Q_m$ are $n \times n$ real symmetric matrices such that 
$$\sum_{k=1}^m \|Q_k\| < 1.$$
Then 
$$\aligned \det^{\quad -{1 \over 2}} \left(I - \sum_{k=1}^m z_k Q_k \right) =&\sum_{k_1, \ldots, k_m \geq 0} {z_1^{k_1} \cdots z_m^{k_m} \over k_1! \cdots k_m!} \\ &\quad \times  {1 \over (2 \pi)^{n/2}} \int_{{\Bbb R}^n} q_1^{k_1}(x) \cdots q_m^{k_m}(x)  e^{-\|x\|^2/2} \ dx, \endaligned \tag3.1.1$$
for all $z_1, \ldots, z_m \in {\Bbb C}$ such that 
$$|z_1|, \ldots, |z_m| \ < \ 1. \tag3.1.2$$
Here we take the principal branch of $\det^{ -{1 \over 2}}$ in the left hand side of (3.1.1), which is equal to $1$ when $z_1=\ldots = z_m=0$. The series in the right hand side converges absolutely and uniformly on compact subsets of the polydisc
(3.1.2).
\endproclaim
\demo{Proof} For $z=\left(z_1, \ldots, z_m\right)$, let 
$$Q_z=I - \sum_{k=1}^m z_k Q_k$$ and let 
$$q_z(x)={1 \over 2} \langle Q x, x \rangle =  {\|x\|^2 \over 2} -\sum_{k=1}^m z_k q_k(x).$$
If $z_1, \ldots, z_m$ are real and satisfy (3.1.2), then $q_z: {\Bbb R}^n \longrightarrow {\Bbb R}$ is a positive definite quadratic 
form, and, as is well known, 
$${1 \over (2 \pi)^{n/2}} \int_{{\Bbb R}^n} e^{-q_z(x)} \ dx ={1 \over \sqrt{\det Q_z}}.$$
Since both sides of the above identity are analytic in the domain (3.1.2), we obtain
$$\det^{\quad -{1 \over 2}} \left(I - \sum_{k=1}^m z_k Q_k \right) ={1 \over (2 \pi)^{n/2}} \int_{{\Bbb R}^n} 
\exp\left\{ -{\|x\|^2 \over 2} + \sum_{k=1}^m z_k q_k(x) \right\} \ dx.$$
Expanding the integral in the right hand side into the series in $z_1, \ldots, z_m$, we complete the proof.
{\hfill \hfill \hfill} \qed
\enddemo

Next, we extract the integral (1.1.2) from the generating function of Lemma 3.1.
Let 
$${\Bbb S}^1 =\left\{z \in {\Bbb C}:\ |z| =1\right\}$$
the the unit circle and let
$${\Bbb T}^m = \underbrace{{\Bbb S}^1 \times \cdots \times {\Bbb S}^1}_{\text{$m$ times}}$$
be the $m$-dimensional torus endowed with the uniform (Haar) probability measure 
$$\mu= \mu_1 \times \cdots \times\mu_m,$$
where $\mu_k$ is the uniform probability measure on the $k$-th copy of ${\Bbb S}^1$. 
If $s \in {\Bbb Z}^m$, $s=(s_1, \ldots, s_m)$, then for the Laurent monomial 
$$\zz^s = z_1^{s_1} \cdots z_m^{s_m},$$ 
we have 
$$\int_{{\Bbb T}^m} \zz^s \ d \mu = \cases 1 & \text{if\ } s=0 \\ 0 &\text{if\ } s \ne 0. \endcases$$
\proclaim{(3.2) Lemma} Let $q_1, \ldots, q_m: {\Bbb R}^n \longrightarrow {\Bbb R}$ be quadratic forms,
$$q_k(x)={1 \over 2} \langle Q_k x, x \rangle \quad \text{for} \quad k=1, \ldots, m,$$
where $Q_1, \ldots, Q_m$ are $n \times n$ real symmetric matrices such that
$$\sum_{k=1}^m\left\| Q_k\right\| \ < \ 1.$$

Then for every $\omega \in {\Bbb C}$ such that $|\omega| < 1$ we have 
$$\split &{1 \over (2 \pi)^{n/2}} \int_{{\Bbb R}^n} \left(1 + \omega q_1(x)\right) \cdots \left(1 + \omega q_m(x)\right) e^{-\|x\|^2/2} \ dx 
\\&=\int_{{\Bbb T}^m} \prod_{k=1}^m \left(1 + \omega z_k^{-1}\right)  
\prod_{(k_1, \ldots, k_s)} \left(1 +  {1 \over 2} \tr \left(Q_{k_1} \cdots Q_{k_s}\right) z_{k_1} \cdots z_{k_s}\right)  d \mu, \endsplit$$
where the second product is taken over all non-empty ordered tuples $(k_1, \ldots, k_s)$ of distinct indices from $\{1, \ldots, m\}$.
\endproclaim
\demo{Proof} From Lemma 3.1, we have
$$\split &{1 \over (2 \pi)^{n/2}} \int_{{\Bbb R}^n} \left(1 + \omega q_1(x)\right) \cdots \left(1 + \omega q_m(x)\right) e^{-\|x\|^2/2} \ dx \\&=
\int_{{\Bbb T}^m} \prod_{k=1}^m \left(1 + \omega z_k^{-1}\right) \det^{\quad -{1 \over 2}} \left( I -\sum_{k=1}^m  z_k Q_k \right) \ d \mu. \endsplit \tag3.2.1$$
Next, we write
$$\split &\det^{\quad -{1 \over 2}} \left( I -\sum_{k=1}^m  z_k Q_k \right) = 
\exp\left\{ -{1 \over 2} \ln \det \left(I - \sum_{k=1}^m z_k Q_k \right) \right\}\\&=
\exp\left\{-{1 \over 2} \tr \ln \left(I - \sum_{k=1}^m  z_k Q_k \right) \right\} =
\exp\left\{ {1 \over 2} \sum_{s=1}^{\infty} {1\over s} \tr \left(\sum_{k=1}^m z_k Q_k \right)^s \right\} \\
&=\exp\left\{  \sum_{s=1}^{\infty} {1 \over 2s}  \sum_{1 \leq k_1, \ldots, k_s \leq m}
 \tr \left(Q_{k_1} \cdots Q_{k_s}\right) z_{k_1}  \cdots z_{k_s}  \right\} \\
 &=\prod_{1 \leq k_1, \ldots, k_s \leq m} \exp\left\{  \sum_{1 \leq k_1, \ldots, k_s \leq m} {1 \over 2s} \tr \left(Q_{k_1} \cdots Q_{k_s}\right) z_{k_1} \cdots z_{k_s}\right\}
\endsplit $$
where the series converges absolutely and uniformly on ${\Bbb T}^m$. 

We expand each of the exponential functions into the Taylor series and observe that only square-free monomials in $z_1, \ldots, z_m$ contribute to the integral (3.2.1), from which it follows that 
$$\split &{1 \over (2 \pi)^{n/2}} \int_{{\Bbb R}^n} \left(1 + \omega q_1(x)\right) \cdots \left(1 + \omega q_m(x)\right) e^{-\|x\|^2/2} \ dx \\&=
\int_{{\Bbb T}^m} \prod_{k=1}^m \left(1+ \omega z_k^{-1}\right) \prod_{(k_1, \ldots, k_s)} \left(1 + {1 \over 2s} \tr \left(Q_{k_1} \cdots Q_{k_s}\right) z_{k_1} \cdots z_{k_s}\right) \ d \mu,
\endsplit $$
where the second product is taken over all non-empty ordered tuples of distinct indices $k_1, \ldots, k_s  \in\{1, \ldots, m\}$.
{\hfill \hfill \hfill} \qed
\enddemo

Our next goal is to write the integral in Lemma 3.2 as the value of the independence polynomial of an appropriate (large) graph. 
\subhead (3.3) Independent sets in weighted graphs \endsubhead
Let $G=(V, E)$ be a finite undirected graph with set $V$ of vertices, set $E$ of edges and without loops or multiple edges. A set $S \subset V$ of vertices is called {\it independent},
if no two vertices from $S$ span an edge of $G$. We agree that $S=\emptyset$ is an an independent set.

Let $w: V \longrightarrow {\Bbb C}$ be a function assigning to each vertex a complex {\it weight} $w(v)$. We define the {\it independence polynomial} of $G$ by
$$\ind_G(w)=\sum\Sb S \subset V \\ S \text{\ independent} \endSb \prod_{v \in S} w(v).$$
Hence $\ind_G(w)$ is a multivariate polynomial in complex variables $w(v)$ with constant term 1, corresponding to $S=\emptyset$. 
\proclaim{(3.4) Corollary} Let $q_1, \ldots, q_m: {\Bbb R}^n \longrightarrow {\Bbb R}$ be quadratic forms,
$$q_k={1 \over 2} \langle Q_k x, x \rangle \quad \text{for} \quad k=1, \ldots, m,$$
where $Q_k$ are real symmetric $n \times n$ matrices and let $\omega \in {\Bbb C}$ be a complex number. 

We define a weighted graph $G=(V, E; w)$ as follows.
The vertices of $G$ are all non-empty ordered tuples $\left(k_1, \ldots, k_s\right)$ of indices $k_1, \ldots, k_s \in \{1, \ldots, m\}$ and two vertices span an edge of $G$ if they have at least one common index $k$, in arbitrary positions. We define the weight of the vertex $\left(k_1, \ldots, k_s\right)$ by 
$${\omega^s \over 2s} \tr\left(Q_{k_1} \cdots Q_{k_s}\right).$$
Then 
$${1 \over (2 \pi)^{n/2}} \int_{{\Bbb R}^n} \left(1 + \omega q_1(x)\right) \cdots \left(1 + \omega q_k(x)\right) e^{-\|x\|^2/2} \ dx = \ind_G(w). \tag3.4.1$$
\endproclaim
\demo{Proof} From Lemma 3.2 it follows that (3.4.1) holds provided $|\omega|$ and $\|Q_k\|$ for $k=1, \ldots, m$ are small enough. Since both sides of (3.4.1) are polynomials in 
$Q_1, \ldots, Q_k$ and $\omega$, the proof follows.
{\hfill \hfill \hfill} \qed
\enddemo

The following criterion provides a sufficient condition for $\ind_G(w) \ne 0$ for an arbitrary weighted graph $G$.
The result is known as the Dobrushin criterion and also as the Koteck\'y - Preiss condition for the cluster expansion, see, for example, Chapter 5 of \cite{FV18}.

\proclaim {(3.5) Lemma} Given a graph $G=(V, E)$ and a vertex $v \in V$, we define its neighborhood $\NN_v \subset V$ by
$$\NN_v= \{v\} \cup \{u \in V: \ \{u, v\}  \in E\}.$$
Let $w: V \longrightarrow {\Bbb C}$ be an assignment of complex weights to the vertices of $G$.
Suppose that there is a function $\rho: V \longrightarrow {\Bbb R}_+$ with positive real values such that for every vertex $v\in V$, we have 
$$\sum_{u \in \NN_v} |w(u)| e^{\rho(u)} \ \leq \ \rho(v).$$
Then 
$$\ind_G(w) \ne 0.$$
\endproclaim
\demo{Proof} See, for example, Section 5.2 of \cite{CF16} for a concise exposition. 
{\hfill \hfill \hfill} \qed
\enddemo

Now we are ready to prove Theorem 1.3. 
\subhead (3.6) Proof of Theorem 1.3 \endsubhead Let $Q_1, \ldots, Q_m$ be the matrices of the quadratic forms $q_1, \ldots, q_m$, so that 
$$q_k(x)={1 \over 2} \langle Q_k, x \rangle \quad  \text{and} \quad \|Q_k \| \ \leq \ {2 \over r} \quad \text{for} \quad k =1, \ldots, m. $$
Since each quadratic form $q_k$ depends of at most $r$ variables, we have
$$\rk Q_k \ \leq \ r \quad \text{for} \quad k=1, \ldots, m.$$
In particular,
$$\left| \tr \left(Q_{k_1} \cdots Q_{k_s}\right)\right| \ \leq \ r \left({2 \over r}\right)^{s} ={2^s \over r^{s-1}}. \tag3.6.1$$
Since each quadratic form $q_k$ has a common variable with at most $r$ other forms, we have
\bigskip
\noindent (3.6.2) For every $k$ there are at most $r$ indices $j \ne k$ such that $Q_k Q_j \ne 0$.
\bigskip 
Let $\omega \in {\Bbb C}$ be a complex number satisfying 
$$|\omega| \leq \gamma={1 \over 4} e^{-{1 \over 2}} \approx 0.1516326649. \tag3.6.3$$

Given $Q_1, \ldots, Q_k$ and $\omega$, we construct a weighted graph $G=(V, E; w)$ as in Corollary 3.4. Our goal is to prove that 
$\ind_G(w)\ne 0$, for which we use Lemma 3.5.

We say that the {\it level} of a vertex $v=\left(k_1, \ldots, k_s\right)$ is $s$ for $s=1, \ldots, m$. Thus for the weight of $v$, we have
$$w(v)={\omega^s \over 2s} \tr\left(Q_{k_1} \cdots Q_{k_s}\right).$$

Combining (3.6.1) and (3.6.3), we conclude that for a vertex of level $s$, we have 
$$|w(v)| \ \leq \ {1 \over s 2^{s+1} r^{s-1}} e^{-{1 \over 2} s}. \tag3.6.4$$
We observe that there are at most $sq r^{q-1}$ vertices $u$ of level $q$ with $w(u) \ne 0$ that are neighbors of a given vertex $v$ (for $q=s$, we count $v$ as its own neighbor). Indeed, there are at most $s$ ways to choose a common index $k$, after which there are at most $q$ positions to place $k$ in $u$. By (3.6.2), we conclude that there are at most 
$sq r^{q-1}$ vertices $u \in \NN_v$ of level $q$ with $w(u) \ne 0$. Choosing $\rho(v) =s/2$ for a vertex of level $s$ and using (3.6.4), we conclude that for a vertex $v$ of level $s$, we have 
$$\split \sum_{u \in \NN_v} |w(u)| e^{\rho(u)} \ \leq \ &\sum_{q=1}^m \left({1 \over q 2^{q+1} r^{q-1} }  e^{-{1 \over 2} q}\right) \left(s q r^{q-1}\right) e^{{1\over 2} q}\\ =
&s \sum_{q=1}^m {1 \over 2^{q+1}} \ < \ {s \over 2}=\rho(v),\endsplit$$
and the proof follows by Corollary 3.4. and Lemma 3.5.
{\hfill \hfill \hfill} \qed

\head 4. Approximating the integral \endhead 

The interpolation method is based on the following simple observation.
\proclaim{(4.1) Lemma} Let $p: {\Bbb C} \longrightarrow {\Bbb C}$ be a polynomial, 
$$p(z)=\sum_{s=0}^m c_s z^s,$$
and $\beta > 1$ be a real number such that 
$$p(z) \ne 0 \quad \text{provided} \quad |z| < \beta.$$
Let us choose a branch of $f(z)=\ln p(z)$ for $|z| < \beta$ and let 
$$T_k(z)=f(0) + \sum_{s=1}^k {f^{(s)}(0) \over s!} z^s$$
be the Taylor polynomial of degree $k$ of  $f$ computed at $z=0$. Then 
$$\left| f(1)-T_k(1)\right| \ \leq \ {m \over  (k+1)\beta^k(\beta-1)}.$$
Moreover, the values of $f^{(s)}(0)$ for $s=1, \ldots, k$ can be computed from the coefficients $c_s$ for $s=0, \ldots, k$
in time polynomial in $n$ and $m$.
\endproclaim
\demo{Proof} See, for example, Section 2.2 of \cite{Ba16}. 
{\hfill \hfill \hfill} \qed
\enddemo
As follows from Lemma 4.1, if $\beta >1$ is fixed in advance, to estimate the value of $f(1)$ within additive error $0< \epsilon <1$ (in which case we say that we estimate the the value of $p(1)=e^{f(1)}$ within relative error $\epsilon$), it suffices to compute the coefficients 
$c_s$ with $s=O\left(\ln n - \ln \epsilon\right)$, where the implied constant in the ``$O$" notation depends only on $\beta$. A similar result holds if $p(z) \ne 0$ in an arbitrary, fixed in advance, connected open set $U \subset {\Bbb C}$ such that $\{0, 1\} \subset U$, see Section 2.2 of \cite{Ba16} (in Lemma 4.1, the neighborhood $U$ is the disc of radius $\beta$).

\subhead (4.2) Computing the integrals \endsubhead Let us fix a constant 
$$0 \ < \ \gamma' \ < \ \gamma,$$
where $\gamma$ is the constant of Theorem 1.3 (so one can choose $\gamma'=0.15$). Let
$q_1, \ldots, q_m: {\Bbb R}^n \longrightarrow {\Bbb R}$ be quadratic forms, defined by their matrices $Q_1, \ldots, Q_m$ as in 
(1.1.1), such that each form depends on not more than $r$ variables among $x_1, \ldots, x_n$ and each form has common variables with not more than $r$ other forms. Suppose that 
the bound (1.3.1) holds. We define a univariate polynomial $p: {\Bbb C} \longrightarrow {\Bbb C}$ by 
$$p(z)={1 \over (2 \pi)^{n/2}} \int_{{\Bbb R}^n} \left(1 + z q_1(x)\right) \cdots \left(1 + z q_m(x)\right) 
e^{-\|x\|^2/2} \ dx.$$
Hence $\deg p \leq m$ and by Theorem 1.3 we have 
$$p(z) \ne 0 \quad \text{provided} \quad |z| < \beta \quad \text{where} \quad \beta = {\gamma \over \gamma'} > 1.$$
In view of Lemma 4.1, to approximate 
$$p(1) = {1 \over (2 \pi)^{n/2}} \int_{{\Bbb R}^n} \left(1 + q_1(x)\right) \cdots \left(1 + q_m(x)\right) e^{-\|x\|^2/2} \ dx \tag4.2.1$$
within relative error $0 < \epsilon < 1$, it suffices to compute 
$p(0)=1$ and $p^{(s)}(0)$ for $s=O(\ln m - \ln \epsilon)$, where the implied constant in the ``$O$" notation is absolute.
From Corollary 3.4, we have 
$$\split p^{(s)}(0)=& s! \sum\Sb \left(k_{11}, \ldots, k_{1s_1}\right), \ldots, 
\left(k_{j1}, \ldots, k_{js_j}\right): \\ s_1 + \ldots + s_j =s \endSb 
{1 \over 2s_1} \cdots {1 \over 2s_j} \\
&\qquad \times \tr \left( Q_{k_{11}} \cdots Q_{1s_1}\right) \cdots \tr \left(Q_{k_{j1}} \cdots Q_{k_{j_s}} \right),\endsplit$$
where the sum is taken over all unordered collections of pairwise disjoint ordered tuples $\left(k_{11}, \ldots, k_{1s_1}\right), \ldots, 
\left(k_{j1}, \ldots, k_{js_j}\right)$ of distinct indices $k_{ij}$ from the set \newline $\{1, \ldots, m\}$, with the total number $s$ of chosen indices. A crude upper bound for the number of such collections is $(2m)^s$: writing all the indices $k_{ij}$ as a row, we have at most $2m$ choices for each index $k_{ij}$, including the choice on whether the index remains in the current tuple or starts a new one. Given that $s=O(\ln m -\ln \epsilon)$ and that computing the traces of the products of $n \times n$ matrices can be done in 
 $(ns)^{O(1)}$ time, we obtain an algorithm approximating the integral in quasi-polynomial 
$n^{O(1)} m^{O(\ln m - \ln \epsilon)}$ time.


\Refs
 \widestnumber\key{AAAA}
 
 
\ref\key{B+11}
\by V. Baldoni, N. Berline, J.A. De Loera, M. K\"oppe and M. Vergne
\paper How to integrate a polynomial over a simplex
\jour Mathematics of Computation
\vol 80
\yr 2011
\pages no. 273, 297--325
\endref

 \ref\key{Ba93}
 \by A. Barvinok
 \paper Feasibility testing for systems of real quadratic equations
 \jour Discrete $\&$ Computational Geometry.
 \vol 10
 \yr 1993
 \pages no. 1, 1--13
 \endref

\ref\key{Ba02}
\by A. Barvinok
\book A Course in Convexity
\bookinfo Graduate Studies in Mathematics, 54
\publ  American Mathematical Society
\publaddr Providence, RI
\yr 2002
\endref 
 
\ref\key{Ba16}
\by A. Barvinok
\book Combinatorics and Complexity of Partition Functions
\bookinfo Algorithms and Combinatorics, 30
\publ Springer
\publaddr Cham
\yr 2016 
\endref


\ref\key{B+19}
\by S. Bravyi, D. Gosset, R. Movassagh
\paper Classical algorithms for quantum mean values
\paperinfo preprint {\tt arXiv:1909.11485}
\yr 2019
\endref

\ref\key{CH88}
\by G.M. Crippen and T.F. Havel
\book Distance Geometry and Molecular Conformation
\bookinfo Chemometrics Series, {\bf 15}
\publ Research Studies Press, Ltd., Chichester; John Wiley $\&$ Sons, Inc.
\publaddr New York
\yr 1988
\endref

\ref\key{CF16}
\by P. Csikv\'ari and P. Frenkel
\paper Benjamini-Schramm continuity of root moments of graph polynomials
\jour European Journal of Combinatorics
\vol 52 
\yr 2016
\pages  part B, 302--320
\endref

\ref\key{FV18}
\by S. Friedli and Y. Velenik
\book Statistical Mechanics of Lattice Systems. A concrete mathematical introduction
\publ Cambridge University Press
\publaddr Cambridge
\yr 2018
\endref

 \ref\key{GP05}
 \by D. Grigoriev and D.V. Pasechnik
 \paper Polynomial-time computing over quadratic maps. I. Sampling in real algebraic sets
 \jour Computational Complexity 
 \vol 14 
 \yr 2005
 \pages no. 1, 20--52
 \endref

\ref\key{H+19}
\by A. Harrow, S. Mehraban and M. Soleimanifar
\paper Classical algorithms, correlation decay, and complex zeros of partition functions of quantum many-body systems
\paperinfo preprint {\tt arXiv:1910.09071}
\yr 2019
\endref
 
\ref\key{LV07}
\by L. Lov\'asz and S. Vempala
\paper The geometry of logconcave functions and sampling algorithms
\jour Random Structures $\&$ Algorithms 
\vol 30 
\yr 2007
\pages no. 3, 307--358
\endref

\ref\key{L+14}
\by L. Liberti, C. Lavor, N. Maculan and A. Mucherino
\paper Euclidean distance geometry and applications
\jour SIAM Review 
\vol 56 
\yr 2014
\pages no. 1, 3--69
\endref

\ref\key{Me04}
\by M.L. Mehta
\book Random Matrices. Third edition
\bookinfo Pure and Applied Mathematics (Amsterdam), 142
\publ  Elsevier/Academic Press 
\publaddr Amsterdam
\yr  2004
\endref

\ref\key{YL52}
\by C.N. Yang and T.D. Lee
\paper Statistical theory of equations of state and phase transitions. I. Theory of condensation. 
\jour Physical Review (2) 
\vol 87 
\yr 1952
\pages 404--409
\endref


\endRefs
\enddocument